\RequirePackage{etoolbox}
\csdef{input@path}{{style/}{graphics/}}
\documentclass[MSNbibl,number,citesort,dvips]{arxbj}
\usepackage{graphicx}
\aid{0}
\volume{21}
\issue{1}
\pubyear{2015}
\firstpage{401}
\lastpage{419}
\doi{10.3150/13-BEJ572} % kopijuoti is 'New paper accepted'

\makeatletter

\newcommand{\rrvert}{\vert}
\newcommand{\llvert}{\vert}
\newcommand{\sgn}{\operatorname{sgn}}
\newcommand{\cov}{\operatorname{cov}}
\newcommand{\RR}{\operatorname{RR}}
\newcommand{\SE}{\operatorname{SE}}
\newcommand{\SD}{\operatorname{SD}}
\makeatother

\begin{document}
\begin{frontmatter}

\title{Beyond first-order asymptotics for Cox regression}
\runtitle{Asymptotics for Cox regression}

\begin{aug}
%%%% inicialai - be tarpu
\author[1]{\inits{D.A.}\fnms{Donald A.} \snm{Pierce}\corref{}\thanksref{1}\ead[label=e1]{pierce.don.a@gmail.com}}
\and
\author[2]{\inits{R.}\fnms{Ruggero} \snm{Bellio}\thanksref{2}\ead[label=e2]{ruggero.bellio@uniud.it}}
%%\runauthor{} %% auto
\address[1]{Department of Public Health, 3181 S.W. Sam Jackson Park
Road, Mail code GH153, Portland, OR 97239-3098, USA. \printead{e1}}
\address[2]{Department of Economics and Statistics, via Tomadini, 30/A
IT-33100, Udine, Italy.\\ \printead{e2}}
\end{aug}

% HISTORY:
\received{\smonth{4} \syear{2012}}
\revised{\smonth{8} \syear{2013}}

% ABSTRACT
%
\begin{abstract}
To go beyond standard first-order asymptotics for Cox regression, we
develop parametric bootstrap
and second-order methods. In general, computation of $P$-values beyond
first order requires more
model specification than is required for the likelihood function. It is
problematic to specify a
censoring mechanism to be taken very seriously in detail, and it
appears that conditioning on censoring
is not a viable alternative to that. We circumvent this matter by
employing a \emph{reference censoring model},
matching the extent and timing of observed censoring. Our primary
proposal is a parametric bootstrap method
utilizing this reference censoring model to simulate inferential
repetitions of the experiment. It is shown
that the most important part of improvement on first-order methods --
that pertaining to fitting nuisance
parameters -- is insensitive to the assumed censoring model. This is
supported by numerical comparisons of
our proposal to parametric bootstrap methods based on usual random
censoring models, which are far more
unattractive to implement. As an alternative to our primary proposal,
we provide a second-order method
requiring less computing effort while providing more insight into the
nature of improvement on first-order
methods. However, the parametric bootstrap method is more transparent,
and hence is our primary proposal.
Indications are that first-order partial likelihood methods are usually
adequate in practice, so we are
not advocating routine use of the proposed methods. It is however
useful to see how best to check on
first-order approximations, or improve on them, when this is expressly desired.
\end{abstract}

% KEYWORDS
% visi is mazosios raides ir pagal abecele
%

\end{frontmatter}

%s1 #&#
\section{Introduction}\label{sec1}
Generally, inferences beyond first order require fuller specification
of the probability model than is needed for the likelihood function
and first-order methods; see Cox and Hinkley \cite{cohi}, Section~2.3. For survival data, a rather broad condition on censoring
mechanisms referred to as `independent censoring' is adequate to
allow computation of the likelihood function and from this many forms
of first-order inference. In principle, methods beyond first order
require more precise specification of the censoring mechanism. 
A~problem is that usual random censoring models are seldom intended to be
realistic, while at the same time implementation of them for going
beyond first-order methods can be quite cumbersome.

We propose use of a \emph{reference censoring model} for bootstrap and
related purposes. This matches the general extent and timing of
observed censoring, but differs from the type of models that are
customarily specified -- while not actually required -- for first-order
methods. In particular, this reference model is progressive Type II
censoring; see, for example, Crowley \cite{crow},
Kalbfleisch and Prentice \cite{kabl}, Section~3.2,
Lawless \cite{lawl}, Section~2.2.1.3, where a fixed number
are censored following each failure, with these fixed numbers matching
the analysis dataset. Jiang and Kalbfleisch \cite{jika}, in mutually
independent work, proposed use of this same reference censoring model.
Its primary virtues are in matching the observed pattern of censoring,
and that implementing it for inferential purposes
beyond first order requires only the rank-based summary data that is
sufficient for partial likelihood, which is not the case for customary
random censoring models. We give considerable attention to the possible
effects of discrepancy between some `actual', or customarily
assumed, censoring model and the reference censoring model.

Our primary proposal is a parametric bootstrap method, detailed in
Section~\ref{sec2.2}. For an hypothesized value of the interest parameter $\psi
$, and an associated estimate of the nuisance parameters, one can
generate bootstrap data under this reference censoring model.
These parameters pertain only to the relative risk function, e.g. $\exp
(z_i \theta)$, in contrast to the situation of the following paragraph.
Inference is based on
tail frequencies in bootstrap trials
of likelihood ratios for testing the hypothesis on $\psi$.
The null distribution of $P$-values
would to third order be uniform $[0, 1]$ if the true censoring model
were our reference model.

This is compared to a more direct parametric bootstrap approach
proposed by
Davison and Hinkley \cite{dahi}, page 351, Algorithms 7.2, 7.3, based on
more conventional random censoring models. This requires estimating,
for plug-in bootstrap use,
further nuisance parameters pertaining to
a censoring distribution and baseline hazard, issues foreign to Cox
regression that our proposal avoids. The dependence on an estimate of
the baseline hazard is due to assumed censoring models depending
directly on time, rather than only on the rank-based data summary. We
find that inferences are similar under these two methods, while the
Davison and Hinkley approach is more cumbersome and prone to
difficulties in fitting the models
in bootstrap trials.

We further provide in Section~\ref{sec4} a second-order asymptotic method
providing inferences
similar to our main proposal but
with far less computation, at some loss of transparency, while
providing additional insight into the nature of improvement on
first-order methods. This method also relies on the reference censoring
model. Connections, in general, between this second-order approach and
parametric bootstrapping are considered by Davison, Hinkley and Young \cite
{dahy}. Such second-order methods were employed by Pierce and Peters
\cite{pipe2} and Pierce and Bellio \cite{pibe} to elucidate the
dependence of $P$-values on aspects of the model that are not required
for likelihood methods, in particular censoring models and stopping
rules. In the latter paper, they showed that for fully parametric
settings in survival analysis, a specific model for censoring is not
required for second-order inferences. Their argument does not fully
apply for partial likelihood in Cox regression, but it does apply to a
major part of the improvement on first-order methods, namely that part
pertaining to effects of fitting nuisance parameters.

In Section~\ref{sec3}, we compare to our parametric bootstrap proposal the
similar method of Jiang and Kalbfleisch \cite{jika}.
Their method obtains confidence intervals with far less computational
time than our proposal, largely since ours is for testing an hypothesis
and must be numerically inverted for a confidence interval. However,
some approximations are made in the Jiang and Kalbfleisch method that
result in less accurate and less powerful inferences.

Samuelsen \cite{samu} considers ``exact''
inference in Cox regression, suggesting with considerable reservation a
method based on
``exact''
logistic regression involving, in turn, those at risk
just prior to each failure, and another method related to our main
proposal but only applicable for Type II censoring at the end of the
follow-up. In Section~\ref{sec5}, we consider issues that arise in his former
proposal, also of interest for other reasons, and obtain further
evidence that this exact logistic regression approach is not satisfactory.

By Cox regression, we mean partial likelihood for survival data
(Cox \cite{cox}) with the relative risk not being
necessarily loglinear. The setting of interest involves inferences
about functions of
$\theta$ in hazard functions of form $\nu(t;z_i,\theta)=\nu_0(t)
\RR(z_i,\theta)$ involving covariate vectors $z_i$, where often
$\RR(z_i,\theta)=\exp(z_i \theta)$.
Censoring is assumed, as usual, to be `independent', meaning roughly
that conditionally on the past, and on covariates, the failure and
censoring times are independent. See
Kalbfleisch and Prentice \cite{kabl}, Sections~1.3, 6.2, and
Andersen \textit{et al.} \cite{abgk}, Section III.2. Kalbfleisch and
Prentice summarize this usefully as independent censoring meaning that
``the probability of censoring at each time $t$ depends only on the
covariate $x$ [of the failure time model], the observed pattern of
failures and censoring up to time $t$ in the trial, or on random
processes that are independent of the failure times in the trial''.
In this, covariates can play a primary role. Censoring based on an
indicator prognostic of failure, even indirectly, will violate the
independence unless such a covariate is correctly included in the model.

We are only interested in use of the usual partial likelihood
estimation, but with inferences improving on usual first-order
approximations. Our approach is also applicable when there is
stratification on the baseline hazard. We consider inference about a
scalar function $\psi(\theta)$ of the relative risk parameters,
framed in terms of testing any specified value for this, with
confidence intervals to be taken as $\psi$-values not rejected. Our
methods are based on the signed likelihood ratio statistic
%
%e1 #&#
\begin{equation}
\label{eq1} r_\psi=\sgn(\hat\psi-\psi) \bigl[2 \bigl\{\ell(\hat
\theta)-\ell (\hat\theta_\psi)\bigr\} \bigr]^{1/2} ,
\end{equation}
where $\ell(\theta)$ is the partial log likelihood function, and
$(\hat\theta, \hat\theta_\psi)$ are respectively, the
unconstrained maximum partial likelihood estimator and that constrained
by the hypothesis $\psi(\theta)=\psi$. Under the hypothesis the
limiting distribution of $r_\psi$ is standard normal, so first-order
likelihood-based $P$-values are $\Phi(r_\psi)$ and $1-\Phi(r_\psi
)$, and the aim here is to improve on that approximation. For the
parametric bootstrap approach, we use the bootstrap distribution of
$r_\psi$ within the progressive Type II censoring framework, with
censoring design adapted to the analysis dataset. For higher-order asymptotics,
we follow that approach and then
provide along lines usually employed for fully parametric models
modifications to $r_\psi$ that are closer to standard normal.

Our aim is for practical settings, and for our main results we
intentionally avoid those with such small sample size that the behavior
of the inference is dominated by discreteness of partial likelihood or
by infinite parameter estimates. For practical purposes, inadequacy of
first-order methods results more from effects of fitting nuisance
parameters in the relative risk than from very small sample sizes. We
provide some limited simulations indicating that, aside from situations
with relatively large numbers of nuisance parameters, first-order
methods are reasonably accurate.

%s2 #&#
\section{The general considerations}\label{sec2}
%s2.1 #&#
\subsection{Some issues of a conditional approach}\label{sec2.1}

Write $c=(c_0,c_1,\ldots,c_m)$ for the censoring configuration, that
is, the number of censorings between successive failure times. Our
initial aim was to condition on this, but we found that this
conditioning is not a satisfactory way to `eliminate' the need to
specify a censoring model. In the first place, there needs to be
further specification regarding which individuals are censored.
Further, the marginal distribution of $c$ will generally depend on all
the parameters including the baseline hazard. Although this does not
preclude conditioning on $c$, it raises issues regarding the loss of
power due to conditioning on $c$. Having found that such conditioning
does not readily resolve the need for specifying a censoring model,
we turned to use of a reference censoring model for hypothetical (or
bootstrap) repetitions of the experiment, in particular that of a
progressive Type II censoring model adapted to the observed censoring
configuration. This inability to resolve the issue by only conditioning
complicates showing that the proposal performs well under more general
censoring models.

Our reference censoring model is closely related to the classical
Kalbfleisch and Prentice \cite{kabl1} result that the partial likelihood,
as a function of the data, provides exactly the probability
distribution of the rank-based data sufficient for partial likelihood.
That is, the partial likelihood is an ordinary likelihood for the
rank-based reduction of the data, which is useful for several needs in
this paper, most particularly for validating the second-order
asymptotics in Section~\ref{sec4}. There are some details of the Kalbfleisch and
Prentice result, considered further in
Kalbfleisch and Prentice \cite{kabl}, Section~4.7.1,
that we should briefly summarize here, and we quote or paraphrase some
of their writing.

Suppose that $k$ items labeled $(1),\ldots,(k)$
give rise to the observed failure times $t_{(1)} < \cdots< t_{(k)}$
with corresponding covariates $Z_{(1)},\ldots, Z_{(k)}$, and suppose
further that $c_i$ items with unobserved failure times
$t_{i1},\ldots,t_{ic_i}$
are censored in the interval
$[t_{(i)}, t_{(i+1)})$. The sets of possible rank summaries at issue
can be represented as
%
%e2 #&#
\begin{eqnarray}
\label{rel}
\begin{array}{l@{\quad}l}\mbox{(A)} & t_{(1)}  < \cdots < t_{(k)} ;
\\
\mbox{(B)} & t_{(i)} < t_{i1}, \ldots , t_{ic_i} \qquad (i=0,1,
\ldots,k).
\end{array}
\end{eqnarray}
The relations (B) mean that the composition of all the risk sets is
specified in this rank summary, and (A) means that the identity of
the item that fails is specified. By these quantities being ``specified''
means only that the covariate values are known, the meaning of a `rank-based' summary. It is straightforward to compute the probability
of (B), conditional on (A),
in terms involving the baseline
hazard $\lambda_0(t)$. In multiplying this by the probability of
(A), the baseline hazard cancels. The result is that the probability
of the rank-based event (\ref{rel}) is given by the partial likelihood
as a function of the data, so the partial likelihood is an ordinary
likelihood for the rank-based data. More details on this are in the
journal article of Kalbfleisch and Prentice \cite{kabl1} than in the
textbook Kalbfleisch and Prentice \cite{kabl}, Section~4.7.1.

A primary subtlety in this result is that the values $c_i$,
corresponding to our censoring configuration, cannot be interpreted
simply as part of the data summary, but rather these must be fixed in
the censoring mechanism. In this regard, Kalbfleisch and Prentice \cite
{kabl}, Section~4.7.1, note that for a \emph{general censoring
mechanism} the joint probability of relations (\ref{rel}) would depend
on the censoring mechanism and the baseline hazard $\lambda_0(t)$.
Only for progressive Type II censoring with the $c_i$ fixed in advance
will the probability of the rank-based data summary be given by the
partial likelihood function. It is also implicit in their calculations
of the probability of (\ref{rel}) that the individuals to be censored
following each failure are chosen uniformly at random from the
corresponding risk set. The issues of this paragraph may be somewhat
clarified by the development of Alvarez-Andrade, Balakrishnan and Bordes \cite{alva}, who
use martingale arguments to derive the partial likelihood function as
giving the probability of relations (\ref{rel}).

%s2.2 #&#
\subsection{Specifics of our proposed method}\label{sec2.2}

Our proposal is based on a direct simulation of the data-generating
model, though employing a reference censoring model, with the interest
parameter at its hypothesized value and nuisance parameters at the
associated estimates. This is often referred to as a ``parametric
bootstrap'', since it uses nuisance parameter estimates from the
analysis dataset. The reference censoring model
approach is to fix, according to the analysis dataset, the censoring
configuration $c$ for repetitions of the experiment, and assume uniform
probability distributions over each risk set
$R[t_{(i)}^+]$ for which individuals are censored. This is the same
reference set as the ``weighted permutation'' re-sampling
employed by Jiang and Kalbfleisch \cite{jika}, in terms of a special
martingale filtration based on the analysis dataset. In the words of
Jiang and Kalbfleisch, this inferential reference set ``imitates
the observed history [of the analysis dataset]'', and ``reproduces the aggregate failure and censoring patterns at all
[observed] failure and censoring times''. That they must also
assume distributions over each risk set for which individuals are
censored suggests that a martingale approach will not resolve the
difficulties in bona fide conditioning.

The specifics of our proposal, in algorithmic form, are as follows:
\begin{enumerate}[3.]
\item[1.] Taking the interest parameter $\psi$ at the hypothesized value,
first generate an uncensored sample of failure times using any baseline hazard,
modulated by the relative risk corresponding to the constrained
estimator $\hat\theta_\psi$.
Then reduce the data to ranks, that is, the time-sorted covariates and
failure indicators required for partial likelihood, which renders
immaterial the choice of baseline hazard; a suitable choice is to take
this as constant.
\item[2.] For the censoring configuration $c$ of the analysis dataset,
censor at random the specified number of individuals following each
failure, using a uniform distribution over each risk set. This can be
done in terms of ranks obtained in Step 1. Carrying this out involves
removing from subsequent risk sets those that are censored at each stage.
\item[3.] This provides a dataset for one trial of the simulation; and for
the parametric bootstrap we compute the partial likelihood ratio
statistic $r_\psi$
by fitting both the null hypothesis and unrestricted models.
\item[4.] From this simulation, we approximate the desired $P$-values as
the bootstrap simulation frequency of $r_\psi$-values more extreme
than that for the analysis dataset.
\end{enumerate}
Since the partial likelihood under progressive Type II censoring is an
ordinary likelihood for the rank reduction of the data, standard theory
for the likelihood ratio parametric bootstrap applies
(DiCiccio, Martin and Stern \cite{dims}, Lee and Young \cite{leey}).
In particular, if the actual censoring distribution were of this type,
then the distribution of $P$-values under the hypothesis would, to
third order,
have a uniform $[0, 1]$ distribution. In Section~\ref{sec4}, we provide some
results pertaining to the actual censoring mechanism not being of this form.

%s2.3 #&#
\subsection{Relation to the Jiang and Kalbfleisch proposal}\label{sec2.3}
The resampling of Jiang and Kalbfleisch \cite{jika} is equivalent to the
bootstrap resampling obtained in Steps 1 and 2, though in line with
their aims for a confidence interval they use the unconstrained estimator
$\hat\theta$ in Step 1, rather than an hypothesized value of $\psi$
and the constrained estimator $\hat\theta_\psi$. However, Step 1 is
implemented differently in not generating failure times, but by
selecting which individual in each risk set is to fail in the
re-sampling, with probabilities proportional to the relative risk
corresponding to $\hat\theta$. That the conditional distribution of
which individual fails in the risk set follows that model is the basis
for partial likelihood. In this respect, their implementation has the
advantage of lending itself more readily to time-dependent covariables.

Their proposal is very different in regard to Steps 3 and 4 above. With
the aim of a confidence interval from a single bootstrap result, they utilize
an approximate pivotal quantity $P(\psi, \mathrm{data})$, whose
distribution is approximately the same for $\psi$-values of
statistical interest and for all values of the nuisance parameter.
They employ for this purpose the score statistic, estimating by
parametric bootstrap its distribution when $\theta=\hat\theta$ of
the analysis dataset, and from this computing
a confidence interval for $\psi$ by the usual pivotal method. The
basic idea for this type of parametric bootstrap was proposed by
Hu and Kalbfleisch \cite{huka},
and in the published discussion of that paper was criticized on grounds that
$P(\psi, \mathrm{data})$ may not be suitably pivotal to higher order.
They also employ a further approximation in regard to the nuisance
parameter. If the
`pivotal'$P(\psi, \mathrm{data})$ were taken as the usual score
statistic, then the constrained maximum likelihood estimator of the
nuisance parameter would be required for each bootstrap trial. To an
approximation, this can be avoided by using a first-order Taylor's
approximation at $\hat\psi$ of the analysis dataset. This is
computationally much faster than our Step 3, which involves fitting
both the unconstrained and constrained estimator in each bootstrap trial.

So altogether, their proposal involves two approximations: choice of an
approximate pivotal and the expansion to approximate the constrained
maximum likelihood estimator. These result in much faster
bootstrap calculations than ours, and arrives at a confidence interval
rather than a test of a specified hypothesis for
$\psi$. Our second-order asymptotic method is computationally as fast
as their proposal for a single bootstrap result, but provides only a
$P$-value for a specified hypothesis rather than a confidence interval.
The two approximations made in their proposal result in somewhat less
power, reflected by somewhat wider confidence intervals, than when our
hypothesis test is inverted to obtain a confidence interval. We assess
this comparison for an example in the following section.

%s3 #&#
\section{Numerical investigation}\label{sec3}
We consider some embellishment of an example from Brazzale, Davison and Reid
\cite{brdr}, Section~7.7, to illustrate points made so far. For their
example, they specify some random censoring models, and an issue is to
see how use of our reference censoring model for repetitions of the
experiment performs under these more commonly-used censoring models.
They consider some higher-order asymptotics issues, to which we return
in Section~\ref{sec5}.

Their example, pertaining to Cox regression with sample size $n = 20$,
has two parts corresponding to: (a) a single binary covariate taking
its values in $3{:} 1$ ratio, with 12.5\% random censoring, and (b) five
Gaussian covariates, of which one carries the interest parameter, with
30\% random censoring. We add to these extensions regarding sample size
and covariate number, to investigate particular needs of this paper. In
particular, for the setting (a) involving no nuisance parameters, we
extend also to $n = 12$ and $n = 40$. As in their example,
the response times are unit exponential random variates. Half of the
observations on the larger treatment arm are subject to censoring at
times uniformly distributed on [0, 4], resulting in expected 12.5\%
censoring for the entire sample. For extensions of setting (b)
involving nuisance parameters, they consider 4 nuisance parameters for
$n = 20$, and we extend also to $n = 40$ with either 4 or 9 nuisance
parameters defined similarly. In simulations, the values of these
nuisance parameters are zero. As in their example for (b) all
observations, which are distributed as exponential with mean unity
without regard to covariates, are subject to censoring at times
uniformly distributed on [0, 3.25], resulting in expected 30\%
censoring. For Tables~\ref{tab1} and \ref{tab2}, the hypothesis is
that the interest parameter in the log relative risk is zero.

First, we give some indication of the performance of usual first-order
methods using the likelihood ratio statistic. These results are
presented in Table~\ref{tab1}, where our interpretation is that
first-order likelihood ratio methods are reasonably adequate for the settings
with relatively few
nuisance parameters, but less so for the settings where there are a
large number of these in relation to the number of failures, that is,
the settings of lines 4 and 6 of Table~\ref{tab1}. Thus, the main
conclusion we draw regarding first-order methods is that for purposes
of evaluating methods of this paper, it is most useful to consider
settings with moderate sample sizes and relatively many nuisance
parameters in the relative risk.

%
%t1 #&#
\begin{table}
\tablewidth=\textwidth
\tabcolsep=0pt
\caption{Null distribution of first-order $P$-values (based on 50\,000 samples).
Table entries are empirical tail frequencies as percentages}\label{tab1}
\begin{tabular*}{\textwidth}{@{\extracolsep{\fill}}lllllllll@{}}
\hline
Nominal & $<$1\% & $<$2.5\% & $<$5\% & $<$10\% & $>$10\% & $>$5\% &
$>$2.5\% & $>$1\% \\
\hline
$n=12$ & 1.6 & 3.8 & 6.9 & 12.9 & 11.0 & 6.3 & 3.1 & 1.4 \\
$n=20$ & 1.5 & 3.6& 6.6 &12.4& 10.0& 5.3& 2.7 &1.2 \\
$n=40$ & 1.3 & 3.0& 5.9& 11.2& 10.0& 5.0 &2.6& 1.1 \\
$n=20, 4$ NP & 3.1& 5.9& 9.7& 15.9& 14.5 & 8.8& 5.4& 2.9 \\
$n=40, 4$ NP & 1.7& 3.7& 7.1& 12.9& 12.8& 7.3& 4.0& 2.0 \\
$n=40, 9$ NP & 3.0& 5.6& 9.2& 15.0& 15.3& 9.4& 5.9& 3.1\\
\hline
\end{tabular*}
\end{table}
%
%
%t2 #&#
\begin{table}[b]
\tablewidth=\textwidth
\tabcolsep=0pt
\caption{Null distribution of first-order and bootstrap $P$-values
(based on 50\,000 samples)}\label{tab2}
\begin{tabular*}{\textwidth}{@{\extracolsep{\fill}}lllllllll@{}}
\hline
Nominal & $<$1\% & $<$2.5\% & $<$5\% & $<$10\% & $>$10\% & $>$5\% &
$>$2.5\% & $>$1\% \\
\hline
\multicolumn{9}{c}{$n=20, 4$ nuisance parameters}\\
$r_\psi$ & 3.1& 5.9& 9.7& 15.9& 14.5 & 8.8& 5.4& 2.9 \\
Reference CM & 1.0& 2.6& 5.1& 10.3& 10.1 & 5.1& 2.6& 1.0 \\
Davison \& Hinkley & 1.2& 2.9& 5.6& 10.9& 10.6 & 5.5& 2.9& 1.1 \\[3pt]
\multicolumn{9}{c}{$n=40, 9$ nuisance parameters}\\
$r_\psi$ & 3.0& 5.6& 9.2& 15.0& 15.3& 9.4& 5.9& 3.1\\
Reference CM & 1.1& 2.6& 5.2& 10.2& 10.1& 5.2& 2.5& 1.0\\
Davison \& Hinkley & 1.1& 2.7& 5.4& 10.4& 10.2& 5.3& 2.7& 1.1\\
\hline
\end{tabular*}
\end{table}

For settings with no nuisance parameters, or relatively few, these
results are consonant with those of Johnson \textit{et al.} \cite{john}
that are in terms of Wald-type inferences rather than likelihood ratio.
As usual,
such results for Wald-type inferences are hampered by the lack of
invariance to parametrization.

Table~\ref{tab2} illustrates the performance of our reference
censoring model parametric bootstrap proposal (`reference CM'),
and the random censoring bootstrap proposal of Davison and Hinkley \cite
{dahi} for the setting of rows 4 and 6 of
Table~\ref{tab1}. For the latter, the function \texttt{censboot} of
the \texttt{R} (R Core Team \cite{R})
package \texttt{boot} (Canty and Ripley \cite{canty})
is used with the option \texttt{sim = ``cond''}, meaning that observed
censoring times are used as potential censoring times in bootstrap trials.
In this paper, for all such simulations, that is, parametric
bootstrapping, we use 10\,000 trials.

At the end of Section~\ref{sec2.2}, it was noted that when the true censoring
model is progressive Type II, the distribution of $P$-values is to
third order
uniform $[0,1]$. The indication from our simulations is that the
accuracy maintains for usual random censoring models, and some
theoretical basis for that is given at the end of Section~\ref{sec4}.
A program employing the routines \texttt{censboot} and \texttt{coxph}
failed in more than 1\% of the bootstrap trials for about 3\% of the
simulated datasets. The $P$-values were computed by ignoring those
bootstrap trials. Many of these failures, however, reflected only
infinite parameter estimates with convergent likelihood. More
seriously, the $P$-values for a small fraction of the datasets may have
been erroneous, without failures in fitting, and this is explained in
more detail below. We do not think these problems have serious effect
on results in Table~\ref{tab2}, but they became more serious when we
used those routines for fitting under parameter values alternative to
the hypothesis, as in the following. For this purpose, we employed a
routine different from \texttt{censboot}, as explained below.

\begin{figure}

\includegraphics{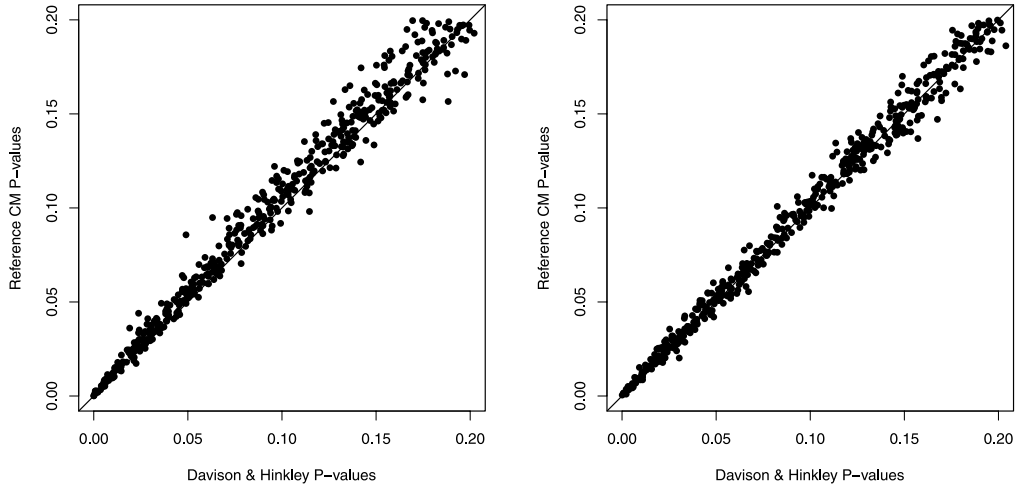}

\caption{Sample-wise comparison of reference censoring model and
Davison \& Hinkley $P$-values, for the setting of $n=20$ with 4
nuisance parameters (left), and $n=40$ with 9 nuisance parameters
(right).}
\label{fig1}
\end{figure}

We indicate in Figure~\ref{fig1} that, sample-by-sample, the $P$-values
from the two proposals in the bottom lines of Table~\ref{tab2} are quite
similar. For each panel of that figure we select, from the calculations
for Table~\ref{tab2}, about 500 analysis samples where the $P$-values
are less than 0.20. Thus the point of our proposal
is far less to improve on the random censoring bootstrap than to obtain
similar results far more easily.
%
%f1 #&#

\begin{figure}

\includegraphics{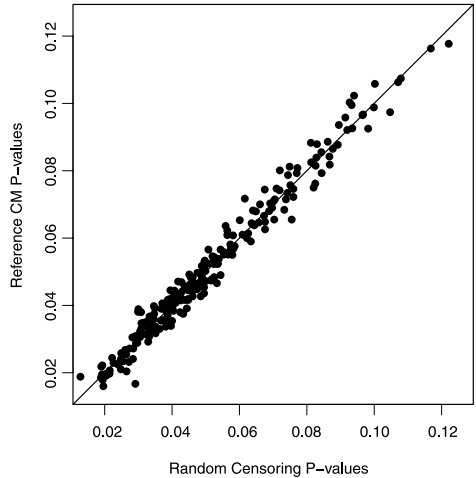}

\caption{Sample-wise comparison of $P$-values from our reference
censoring model proposal and the censoring model based on the
`clinical trial' considerations, in the setting with $n=20$ and 4
nuisance parameters.}
\label{fig2}
\end{figure}

In Figure~\ref{fig2}, we present more limited results of this nature
under the alternative, i.e. when the true value of $\psi$ differs from
the hypothesized value. We employ a different form of censoring model
of some interest, discussed below.
In this figure, we compare $P$-values from our reference censoring
model proposal and those of the bootstrap employing the actual
censoring model, for 240 simulated samples with both methods being used
for each simulation sample.
The true value of $\psi$ is taken as zero and the hypothesis being
tested on each trial is that $\psi$ is equal to the 95\% Wald upper
confidence limit. Again we see that, sample-by-sample, the $P$-values
from the two methods are similar. The value of the alternative is
allowed to vary between simulation trials, since the aim is not to
estimate powers for the two methods, but to show more fundamentally
that the inferences from the two methods are quite similar. For this
aim, it was more effective to be always carrying out tests that have
interesting $P$-values.
%
%f2 #&#

Generally, a major reason for censoring is the end of follow-up,
and this is the mechanism employed for the simulation leading to Figure~\ref{fig2}.
For a prototypical clinical trial setting, patients are
enrolled at random during some period, say the first 2 years.
Follow-up continues for 5 more years, say, at which time all remaining
subjects are censored. For simulation, a constant failure rate value
was chosen to yield 30\% censoring. Analysis time is the interval from
enrollment until failure or censoring. The primary distinction between
this censoring model and the random censoring model employed above is
that, after enrollment is complete, the potential censoring times for
all subjects are known. Thus these can be taken as fixed for the
bootstrap trials, rather than only for the observed censoring times
as in \texttt{censboot}. Nevertheless, implementing this entails use
of an estimate of the baseline hazard, since the censoring model
involves times rather than only ranks.

For the parametric bootstrap using this censoring model, the analog of
the Davison and Hinkley proposal employed so far, we have not used the
routine \texttt{censboot} for results in Figure~\ref{fig2}. This is because
the censoring model differs from that used for Table~\ref{tab2} and
Figure~\ref{fig1}
and it is not required to estimate a censoring distribution. We have
also modified the simulation of failure times for the following reason.
Not surprisingly, with sample sizes as small as $n = 20$ and 30\%
censoring, a method using nonparametric estimation of the
baseline hazard can fail, unless more smoothing is employed than by
\texttt{censboot}. For example, this happens when there are several
censorings following the last failure, so that any nonparametric
estimator of the
baseline hazard
is then undefined after decreasing to a value considerably greater than
zero. The \texttt{censboot} routine places the substantial undefined
mass just after the last failure time, resulting in $P$-values
considerably different than those using the true survival distribution,
or those from our reference censoring model proposal. This problem was
more serious in testing under the alternative for
Figure~\ref{fig2}
than under the null hypothesis of Figure~\ref{fig1}. In the routine
used for Figure~\ref{fig2}, we employed in bootstrap trials somewhat
more smoothing to deal with this, fitting a Weibull distribution to the
incomplete nonparametric estimate of the baseline
survival distribution.

We now offer some comparison to the Jiang and Kalbfleisch \cite{jika}
proposal. For obtaining a confidence interval, their procedure is about
twice as fast than our bootstrap method. Half of this is due to
employing an approximate pivotal so that only one bootstrap run is
required, and the remaining half is due to using an approximate score
to avoid model fitting on bootstrap trials. These approximations were
considered in more detail in Section~\ref{sec2.3}. Our aim is to evaluate the
effect of these approximations for the example above involving $n=20$
with 4 nuisance parameters. The censoring model for data generation was
that described at the beginning of this section, and used for Tables~\ref{tab1} and \ref{tab2}. We also present operating characteristics
for the first-order method based on the normal approximation to the
distribution of
$r_\psi$. Table~\ref{tab3n} presents in column 1 the estimated
probability of coverage for a 95\% lower confidence limit for one of
the relative risk parameters when the true value is zero, based on
simulation of 10\,000 datasets for the Jiang and Kalbfleisch method, and
taking first-order and our bootstrap values from Tables~\ref{tab1} and
\ref{tab2}. In order to evaluate what corresponds to the `length'
of the one-sided confidence interval, we present in column 2 the
probability of the interval covering the false value $-0.50$, based on
$2000$ datasets. This particular false value was chosen to be roughly
in the center of the distribution of lower confidence limits. The
number of trials is smaller than for other purposes here, since the
computations are extensive, so standard errors are reported.
%
%t3 #&#
\begin{table}
\tablewidth=\textwidth
\tabcolsep=0pt
\caption{Coverage probability for 3 methods.
First column: based on 50\,000 samples for $r_\psi$ and reference
censoring model bootstrap, and 10\,000 samples for the Jiang \&
Kalbfleisch method.
Second column: based on 2000 samples}\label{tab3n}
\begin{tabular*}{\textwidth}{@{\extracolsep{\fill}}lll@{}}
\hline
Method & Coverage (\%) true value 0 & Coverage (\%) of $-0.50$ \\
\hline
$r_\psi$ & 88.7 & 58 $\pm$ 1 \\
Jiang \& Kalbfleisch & 93.6 & 64 $\pm$ 1 \\
Reference CM & 94.8 & 42 $\pm$ 1 \\
\hline
\end{tabular*}
\end{table}

The coverage probability in column 1 for the Jiang and Kalbfleisch
method, for which the standard error is 0.25\%,
could be considered as adequate, but as anticipated the method is
slightly inferior compared to our more computationally intensive method
that does not employ an approximate pivotal quantity.
It should be borne in mind that this example was chosen specifically
to place considerable stress on first-order methods, and hence to
challenge methodology for improvement. In their paper, they consider
only examples with $n = 50$ or 60, with either no nuisance parameters or
only one.
In the final column of Table~\ref{tab3n} it is seen that as expected,
due to the approximations resulting in far less computational time,
their confidence intervals are somewhat longer than ours, that is the
probability of covering the false value $\psi=-0.5$ is about 50\% larger.

%s4 #&#
\section{Higher-order asymptotics}\label{sec4}
For the reference censoring model approach, or use of any specified
censoring model, $P$-values that we have above approximated using the
parametric bootstrap can also be approximated using the type of
higher-order asymptotics developed by Barndorff-Nielsen \cite
{barn1,barn2}, Barndorff-Nielsen and Cox \cite{baco} and others. A comprehensive textbook
treatment is given in Severini \cite{seve}.
The basic issues regarding censoring models are the same as for the
parametric bootstrap. For our progressive Type II censoring model, the
higher-order asymptotics are valid for Cox regression partial
likelihood since this corresponds to ordinary likelihood for the
rank-based data. Extension to other censoring models is considered
later in this section, and further in the Discussion.

The higher-order asymptotic methods consist of modifying the directed
likelihood ratio (\ref{eq1}) to have more nearly a
standard normal distribution. This is given by $r^*_\psi$ of (\ref
{eq2}) that is used as in (\ref{rs}) below. To implement this for our
setting involves simulation to approximate certain likelihood
covariances. This differs from the parametric bootstrap in that no
model fitting is involved in the simulation, which can be important
when there are convergence difficulties in fitting. Further, the
required number of simulation trials is smaller than for the bootstrap
since this is for estimating covariances rather than tail probabilities
more directly.
This form of asymptotics has been found in many investigations to be
remarkably accurate, and for our setting this is borne out in numerical
results of Table~\ref{tab3}. It also provides useful insights by
separating the effects of fitting nuisance parameters from those due to
limited adjusted information.

The following results are presented for the more general setting as
developed in Barndorff-Nielsen and Cox \cite{baco}, Section~6.6,
rather than attempting to specialize to our partial likelihood setting.
It should be noted that the formulas are far less to provide recipes
for calculation, than to indicate primary concepts of the methodology.
In particular,
for many of the formulas to make sense, the data must first be
represented to second order as $(\hat\theta,a)$, the maximum
likelihood estimator along with a suitable approximate ancillary. For
this reason, the formulas found little practical use until the
development by Skovgaard \cite{skov1,skov2} of approximations in terms
in terms of the approximate ancillary based on $\hat\imath^{-1} \hat
\jmath$ arising in the formulas to follow.
Fraser, Reid and Wu \cite{fra99} developed a different approach to the
approximation, less general and perhaps difficult to apply to the
present setting.

The Barndorff-Nielsen adjustment (Barndorff-Nielsen and Cox \cite{baco}, Section~6.6), has general form
%
%e3 #&#
\begin{equation}
\label{eq2} r^*_\psi=r_\psi+ r_\psi^{-1}
\log(u_\psi/r_\psi)=r_\psi+\mathrm
{NP}_\psi+ \mathrm{INF}_\psi ,
\end{equation}
where $u_\psi$ involves derivatives of the likelihood function with
respect to parameter estimates, holding fixed a notional ancillary.
This means a suitable complement to the maximum likelihood estimator to
provide a second-order sufficient statistic; see Barndorff-Nielsen and Cox
\cite{baco}, Section~2.5.
We will often suppress the subscript $\psi$, in phrases similar to ``the $r^*$ method''.
Similarly, we will sometimes suppress the subscripts on $\mathrm
{NP}_\psi$ and $\mathrm{INF}_\psi$. The nuisance parameter
adjustment $\mathrm{NP}_\psi$ corresponds to the modified profile
likelihood $L_{\mathrm{MP}}(\psi)$, which was developed in theory not directly
related to $r^*$, and applies to vector parameters $\psi$ as well; see
Barndorff-Nielsen \cite{barn0} and
Barndorff-Nielsen and Cox \cite{baco}, Section~8.2. For scalar parameter
$\psi$, the exact relation is $L_{\mathrm{MP}}(\psi) = \exp(-r_\psi \mathrm
{NP}_\psi) L_P(\psi)$.
The information adjustment $\mathrm{INF}_\psi$
allows for limited adjusted information for $\psi$, more specifically
pertaining to the skewness of the score statistic for
inference about $\psi$; see Pierce and Peters \cite{pipe}.

The quantities in (\ref{eq2}) are specifically
\[
\mathrm{NP}_\psi=r_\psi^{-1} \log(C_\psi)
, \qquad \mathrm{INF}_\psi= r_\psi^{-1} \log(\stackrel{
\times} {u}_\psi /r_\psi) ,
\]
with\vspace*{-2pt}
\begin{eqnarray*}
\stackrel{\times} {u}_\psi& =& \hat{\jmath}^{-1/2}_{\psi\psi|\nu
}
\biggl[ \frac{\partial\{ \ell_P(\psi)- \ell_P(\hat\psi)\}}{\partial\hat
\psi} \biggr] ,
\\
C_\psi& =& \biggl\llvert \frac{\partial^2 \ell(\psi,\hat\nu_\psi)}{
\partial\hat\nu \,\partial\nu} \biggr\rrvert / \{ | \tilde
\jmath_{\nu\nu}| |\hat\jmath_{\nu\nu}| \}^{1/2} .
\end{eqnarray*}
We note that
\begin{eqnarray*}
\frac{\partial\{ \ell_P(\psi)- \ell_P(\hat\psi)\}}{\partial\hat
\psi}  = \frac{\partial\ell(\hat\theta_\psi)}{\partial\hat\psi}-
\frac{\partial\ell(\hat\theta)}{\partial\hat\psi} -\frac{\partial^2 \ell(\hat\theta_\psi)}{\partial\hat\psi\,
\partial\nu}
\biggl\{\frac{\partial^2 \ell(\hat\theta_\psi)}{\partial\hat\nu
\,\partial\nu} \biggr\}^{-1} \biggl\{ \frac{\partial\ell(\hat\theta_\psi
)}{\partial\hat\nu}-
\frac{\partial\ell(\hat\theta)}{\partial
\hat\nu} \biggr\} ;
\end{eqnarray*}
see
Barndorff-Nielsen and Cox \cite{baco}, equation 6.106.
Here $\theta$ is represented as $(\psi,\nu)$, where
$\nu$ is any version of the nuisance parameter, and
$\tilde\jmath_{\nu\nu}$
and $\hat\jmath_{\nu\nu}$ are the observed information matrices for
$\nu$ evaluated at the constrained
and unconstrained maximum likelihood estimates,
$\hat{\jmath}^{-1/2}_{\psi\psi|\nu}$\vspace*{2pt}
is the observed adjusted information for $\psi$ as defined
in Section~9.3(iii) of Cox and Hinkley \cite{cohi}. The partial derivatives
in the above expressions
are generally referred to as sample-space derivatives, where for instance
\[
\frac{\partial\ell(\hat\theta_\psi)}{\partial\hat\psi} = \frac
{\partial\ell(\theta;\hat\psi,\hat\nu,a)}{\partial\hat\psi
} \biggl|_{\theta=\hat\theta_\psi} .
\]

Skovgaard \cite{skov1,skov2}
derived second-order approximations to these quantities that avoid
difficult partial derivatives.
These results
were further developed by Severini \cite{seveemp},\vadjust{\goodbreak}
without the reliance made by Skovgaard on curved
exponential families.
We first express these in rather casual but useful notation as
\begin{eqnarray*}
\frac{\partial\ell(\hat\theta_\psi)}{\partial\hat\theta\,
\partial\theta} & \doteq& \cov_{\hat\theta}\bigl\{ U(\hat\theta), U(\hat
\theta_\psi)\bigr\} \hat\imath^{-1} \hat\jmath ,
\\
\frac{\partial\ell(\hat\theta_\psi)}{\partial\hat\theta}-\frac
{\partial\ell(\hat\theta)}{\partial\hat\theta} &\doteq& \cov _{\hat\theta}\bigl\{ U(
\hat\theta), \Delta(\hat\ell)\bigr\} \hat\imath^{-1} \hat\jmath ,
\end{eqnarray*}
where $U$ is the score $\partial\ell/\partial\theta$, and $\Delta$
is a log likelihood difference defined below.
Here $\hat\jmath$\vspace*{2pt} and $\hat\imath$ are the observed and expected
information matrices for the full parameter, evaluated at
$\hat\theta$. More precisely
%
%e4 #&#
\begin{eqnarray}
\label{eq3} \cov_{\hat\theta}\bigl\{ U(\hat\theta), U(\hat
\theta_\psi)\bigr\} & =& \cov_{\theta_1}\bigl\{ U(
\theta_1), U(\theta_2)\bigr\} ,\nonumber
\\[-8pt]\\[-8pt]
\cov_{\hat\theta}\bigl\{ U(\hat\theta), \Delta(\hat\ell)\bigr\} &= & \cov
_{\theta_1}\bigl\{ U(\theta_1), \ell(\theta_2)-\ell(
\theta_1)\bigr\} ,
\nonumber
\end{eqnarray}
where following computation of the expectations, $\theta_1$ and
$\theta_2$
are respectively, evaluated at the unconstrained and constrained
maximum likelihood estimators for the analysis dataset $\hat\theta$
and $\hat\theta_\psi$.
The covariances in (\ref{eq3}) can be approximated by simulation.
A key aspect of the Skovgaard
result is that the quantity $\hat\imath^{-1} \hat\jmath$ captures
adequately the ancillary information. Quantities above differ from
usual information calculations through involving covariances of scores
at two different parameter values.

For either $r^*_\psi$ or the Skovgaard approximation to it, the sense
of the final approximation can be usefully expressed as follows. For
this, we write
$Y$ as a random dataset,
as $y$ the observed value, and indicate explicitly the dependence of
$r_\psi$ and $r^*_\psi$ on the data. Then we have that
%
%e5 #&#
\begin{equation}
\label{rs} P \bigl\{r_\psi(Y) < r_\psi(y); \psi(\theta)=
\psi \bigr\}= \Phi \bigl\{ r^*_\psi(y) \bigr\} \bigl[ 1+
\mathrm{O}_p \bigl\{ (\hat\psi-\psi) n^{-1/2} \bigr\} \bigr] ,
\end{equation}
which to this order does not depend on the nuisance parameter in the
relative risk. For the present context $n$
is best thought of as the number of failures. Thus, for deviations
$\hat\psi-\psi=\mathrm{O}_p(n^{-1/2})$ the relative error
is $\mathrm{O}(n^{-1})$
and for large deviations the error is $\mathrm{O}(n^{-1/2})$.
The large deviation property is important in ruling out second-order
approximations that suffer from being overly local. The above bound
also holds when conditioning on suitable ancillary statistics,
in particular that based on $\hat\imath^{-1} \hat\jmath$; see
Skovgaard \cite{skov0}.
The justification for (\ref{rs}) is a combination of the result (7.5)
of Severini \cite{seve} and results in Section~4 of
Skovgaard \cite{skov1}.

To approximate the covariances (\ref{eq3})
in our setting one may carry out Steps 1 and 2 of the algorithm given
for the parametric bootstrap, but keeping track only of the required
likelihood quantities so that their covariances can be computed
following the simulation. Steps 3 and 4 are not required, and no model
fitting is involved in this simulation, which avoids not only computing
effort but possible convergence difficulties. Note that in this
simulation the
resampling is done under the unconstrained estimator $\hat\theta$. As
noted, estimation of these covariances requires far fewer bootstrap
trials than for the parametric bootstrap. For accurate results, it is
important that $\hat\imath$ be approximated from the same simulation
samples as the other covariances, even when a formula for it is known.
This is because, as laid out in Sections~2 and 3 of Severini \cite{seveemp}, when $\hat\imath$ is approximated in this
manner the leading terms of what is being approximated,
and the covariance-based approximations, are identical.

We illustrate these methods by continuing with the example of Table~\ref{tab2}.
Recall that the true censoring model for the analysis datasets differs
from our reference progressive Type II model used
for the second-order and bootstrap calculations, and were described
early in Section~\ref{sec3}. Table~\ref{tab3}
indicates that the second-order and parametric bootstrap $P$-values are
quite similar. The results in the last line of that table,
though probably adequate for practice, contain some entries differing
by 6--9 simulation standard errors from nominal values.
We discuss that in the final section.
%
%t4 #&#
\begin{table}
\tablewidth=\textwidth
\tabcolsep=0pt
\caption{Null distribution of $P$-values for the higher-order method based on
$r^*_\psi$ (based on 50\,000 samples)}\label{tab3}
\begin{tabular*}{\textwidth}{@{\extracolsep{\fill}}lllllllll@{}}
\hline
Nominal $\%$ tail & $<$1\% & $<$2.5\% & $<$5\% & $<$10\% & $>$10\% &
$>$5\% & $>$2.5\% & $>$1\% \\
\hline
\multicolumn{9}{c}{$n=20, 4$ nuisance parameters}\\
$r_\psi$ & 3.1& 5.9& 9.7& 15.9& 14.5 & 8.8& 5.4& 2.9 \\
Bootstrap (reference CM) & 1.0& 2.6& 5.1& 10.3& 10.1 & 5.1& 2.6& 1.0 \\
$r^*_\psi$ (reference CM) & 1.0& 2.5& 4.9& 9.7& 9.5& 4.8& 2.5& 1.0\\[3pt]
\multicolumn{9}{c}{$n=40, 9$ nuisance parameters}\\
$r_\psi$ & 3.0& 5.6& 9.2& 15.0& 15.3& 9.4& 5.9& 3.1\\
Bootstrap (reference CM) & 1.1& 2.6& 5.2& 10.2& 10.1& 5.2& 2.5& 1.0\\
$r^*_\psi$ (reference CM) & 0.8& 2.1& 4.3& 9.0& 8.9& 4.3& 2.1& 0.8\\
\hline
\end{tabular*}
\end{table}
%
%f3 #&#
\begin{figure}[b]

\includegraphics{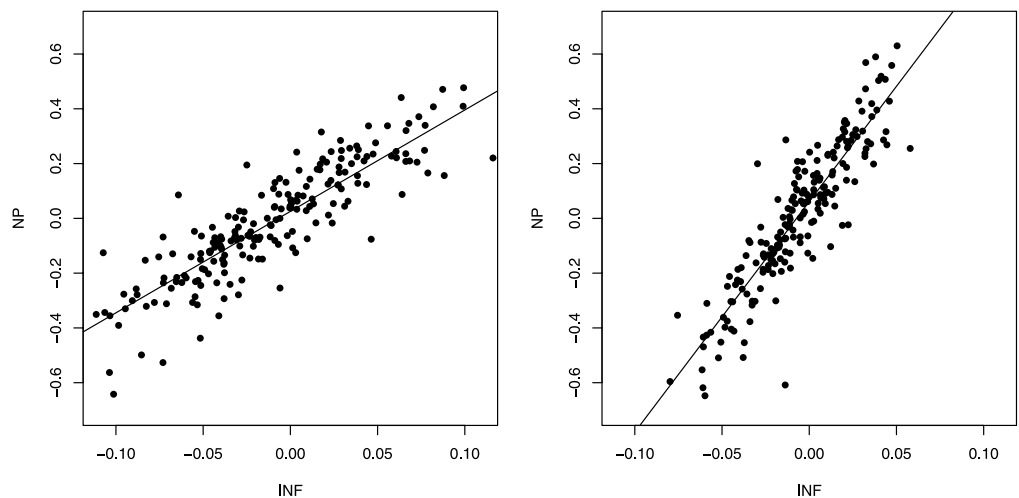}

\caption{Values of $\mathrm{NP}_\psi$ and $\mathrm{INF}_\psi$ for
$200$ trials
for $n=20$ with 4 nuisance parameters and $n=40$ with 9 nuisance
parameters, and for testing the hypothesis that $\psi$ is equal to the
95\% lower Wald confidence limit.
The slopes of the regression lines are 3.7 and 8.4.}
\label{fig3}
\end{figure}

With this method the nature of the improvements on first-order
inference can be characterized more clearly than for the parametric bootstrap.
For each of the examples of Table~\ref{tab2}, we describe in Figure~\ref{fig3} the nature of the
$\mathrm{NP}$ and $\mathrm{INF}$ adjustments for 200 simulated
samples. This simulation was carried out under the strategy used for
Figure~\ref{fig2}; that is, to focus on interesting $P$-values we take
the hypothesis being tested for each sample as the lower 95\% Wald
confidence\vspace*{1pt} limit
$\hat\psi-1.645 \SE(\hat\psi)$. This means that the value of $r_\psi
$ varies only modestly around the value $-1.645$.

The $\mathrm{NP}$ adjustments are in magnitude about 4 and 8 times
larger than the $\mathrm{INF}$ adjustments, for the left and right
panels of Figure~\ref{fig3}, as indicated by the slopes of the
regression lines. This reflects an important point that is not
restricted to the setting of this paper. Unless the sample size is
quite small the
$\mathrm{INF}$ adjustment is usually minor, but when there are many
nuisance parameters, relative to the sample size, the $\mathrm{NP}$
adjustment can be substantial even for moderately large samples.
Whether or not this is the case depends on the particular setting, and
is often difficult to predict without making the computation.
Consideration of this was first given in Section~3 of Pierce and Peters \cite{pipe} and has been further described by several others:
for example, Brazzale, Davison and Reid \cite{brdr} and Sartori
\cite{sart}, Section~6.

Pierce and Bellio \cite{pibe} employed these second-order asymptotics
to investigate the inferential effects of model specifications that do
not affect the likelihood function. They found that for ordinary
likelihood, inferences are usually to second order unaffected by choice
of censoring models. On the other hand, they found that the effect of
stopping rules is to second order carried only by the $\mathrm{INF}$
adjustment, with the $\mathrm{NP}$ adjustment unaffected. Their
argument for censoring models in ordinary likelihood does not fully
apply to partial likelihood for survival data, since it depends on the
observations being stochastically independent. However the argument for
the $\mathrm{NP}$ part of the adjustment only requires that the
contributions to the score are uncorrelated, which is the case for
partial likelihood. See, for example, Cox \cite{cox75}
but note that this is an essential aspect of martingale theory. Many of
the standard martingale developments need substantial modification for
partial likelihood, and this is treated by Fleming and Harrington \cite
{flha}. The point of this argument is that in terms of the second-order
asymptotics, it is not necessary that the true censoring model be the
same as our progressive Type II reference censoring model.

%s5 #&#
\section{Considering risk sets as fixed}\label{sec5}
There is an inferential frame of reference having the same likelihood
function as the Cox regression partial likelihood, which, though not
conditional, has been and remains of considerable interest. This
consists of considering all the risk sets, i.e. those at risk just
before each failure, as fixed in hypothetical repetitions of the
experiment. In particular, this is the basis for the ``logistic
exact'' inference
for Cox regression proposed by Samuelsen \cite{samu},
although with some reservations. The point of this section is to show
that, although a definition of $r^*$ is temptingly very simple for this
frame of reference, it is not suitable for use for survival data, as
the adjustments $r^*-r$ are too small to be of the value we have seen
throughout this paper.

This frame of reference leads to a multinomial probability model
defined on each risk set, with the probabilities of which individual
fails proportional to the relative risk. This seems to be what
Cox \cite{cox} originally had in mind, although it was
quickly realized that this is not a conditional frame of reference.
Indeed, if there were no censoring fixing all the risk sets allows no
data variation for purposes of partial likelihood. Unless censoring is
quite heavy, the risk sets remain to be substantially determined by
failures. These difficulties turn out to have a large bearing on the
higher-order asymptotics. They were resolved by the martingale approach
to Cox regression, in which the risk sets are fixed successively but
not simultaneously.

We have seen throughout this paper that a useful definition of $r^*$
for survival data is possible, and in extreme situations the
improvement over $r$ is substantial. This does not happen with the
definition of
$r^*$ appropriate for the fixed risk set frame of reference. Thus,
although the product multinomial formulation is a probability model for
some experiment, it is unsuitable for going beyond first order in
analysis of survival data.

In Table~\ref{tab4} we show, for the two examples of Table~\ref{tab3}, results of using $r^*$ for the fixed risk set frame of
reference, along with those of the $r^*$ method proposed in Section~\ref{sec4},
labeled ``reference CM'' for the progressive Type II
reference censoring model employed. The fixed-risk-set $r^*$ performs
very little better than the first order $r$, even though the proposal
of Section~\ref{sec4} results in substantial improvement.
%
%t5 #&#
\begin{table}
\tablewidth=\textwidth
\tabcolsep=0pt
\caption{Null distribution of $P$-values for the $r^*$ method with fixed risk
sets (based on 50\,000 samples)}\label{tab4}
\begin{tabular*}{\textwidth}{@{\extracolsep{\fill}}lllllllll@{}}
\hline
Nominal $\%$ tail & $<$1\% & $<$2.5\% & $<$5\% & $<$10\% & $>$10\% &
$>$5\% & $>$2.5\% & $>$1\% \\
\hline
\multicolumn{9}{c}{$n=20, 4$ nuisance parameters}\\
$r_\psi$ & 3.1& 5.9& 9.7& 15.9& 14.5 & 8.8& 5.4& 2.9 \\[2pt]
$r^*_\psi$ (fixed risk set) & 2.6& 5.3& 9.0& 15.2& 13.5& 7.8& 4.7&
2.4\\[4pt]
$r^*_\psi$ (reference CM) & 1.0& 2.5& 4.9& 9.7& 9.5& 4.8& 2.5& 1.0\\[3pt]
\multicolumn{9}{c}{$n=40, 9$ nuisance parameters}\\
$r_\psi$ & 3.0& 5.6& 9.2& 15.0& 15.3& 9.4& 5.9& 3.1\\[2pt]
$r^*_\psi$ (fixed risk set) & 2.5& 5.0& 8.3& 14.0& 15.5& 9.4& 5.9&
3.0\\[4pt]
$r^*_\psi$ (reference CM) & 0.8& 2.1& 4.3& 9.0& 8.9& 4.3& 2.1& 0.8\\
\hline
\end{tabular*}
\end{table}

We now show that the $r^*$ method for the fixed risk-set frame of
reference, even though not usually suitable for Cox regression, is very
simple in form. The resulting product multinomial setting falls in the
simple framework for the higher-order asymptotics laid out by
Pierce and Peters \cite{pipe}, that is, full rank exponential
families in terms of canonical parameters. In this case, we have that
%
%e6 #&#
\begin{eqnarray}
\label{eq4} \mathrm{INF}_\psi& =& r_\psi^{-1}
\log(w_\psi/r_\psi) ,\nonumber
\\[-8pt]\\[-8pt]
\mathrm{NP}_\psi&=& r_\psi^{-1} \log\bigl(
\rho_\psi^{1/2}\bigr) ,
\nonumber
\end{eqnarray}
where $w_\psi$ is the Wald statistic for the canonical multinomial
interest parameter and $\rho_\psi$
is the determinant ratio of canonical nuisance parameter information matrices,
namely
\[
\rho_\psi=\frac{ |j_{\nu\nu}(\hat\theta)|}{|j_{\nu\nu}(\hat
\theta_\psi)|} .
\]
When the relative risk is loglinear as indicated initially in Section~\ref{sec1}, the parameters $\psi$
and $\nu$ are coordinates of $\theta$. Thus the ingredients for
computing are readily available from the usual partial likelihood
fitting. The adjustments (\ref{eq4}) are those considered in
Brazzale, Davison and Reid \cite{brdr}.

Considering the risk sets as fixed can be useful for a parametric
bootstrap when all the risk sets are very large, and their composition
is mainly determined by censoring or competing risks. That is, when the
events under study are rare. The difficulty with considering the risk
sets as fixed is that they are in part determined by failures. This
aspect becomes negligible, though, in settings as just described. Such
settings commonly arise in epidemiology, but less frequently in
clinical trials. What transpires is that the adjustments in (\ref
{eq4}) seem always quite small, which is reflecting statistical issues
for these types of settings.

%s6 #&#
\section{Discussion}\label{sec6}
Any method for going beyond first order for Cox regression must involve
some form of specification of the censoring model. Usual random
censoring models are seldom intended to be realistic, but only to
specify a concrete model compatible with independent censoring. But
even if such a model were realistic, utilizing it for inference beyond
first order would usually involve unattractive estimation of censoring
distributions and the baseline hazard. It appears to us that the only
highly tractable approach is to utilize some kind of reference
censoring model, chosen to match primary aspects of the observed
censoring. We have proposed such a reference censoring model, and shown
that to second order the dominant aspect of adjustment to first-order
methods is independent of the choice of reference censoring model. Our
conclusions in these matters are compatible with those of Jiang and Kalbfleisch \cite{jika}, based on mutually independent work.

We are not advocating routine use of the methods of this paper, as
ordinarily first-order methods seem adequate for practical purposes.
When one does desire some confirmation of this, the second-order
methods of Section~\ref{sec4} involve far less computation than the parametric
bootstrap method of Section~\ref{sec3}, and the results of these two approaches
are similar. The parametric bootstrap method is considerably more
transparent, which we consider important. It appears that, not
surprisingly, in extreme situations the parametric bootstrap results
are slightly more accurate than the second-order ones. This is because
the second-order methods rely on the distribution of $r_\psi$ being
not terribly far from standard normal, which is nearly always the case
in practice. One could probably construct extreme examples where the
improvement of the parametric bootstrap over the second-order method is
greater than we have seen in this paper.

We would not want to leave readers with the impression that the
higher-order asymptotics of Section~\ref{sec4} applies to Cox regression only
for our reference progressive Type II censoring model. Mykland and Ye
(U. Chicago Statistics Dept. TR 332, 1995) showed that under usual
conditions of independent censoring, Bartlett identities of all orders
are satisfied for Cox regression partial likelihood. Mykland
\cite{mykl} showed that this result is adequate to establish that
adjustments of affine form $r^\dagger_\psi=\{ r_\psi- E(r_\psi) \}
/ \SD(r_\psi)$
are standard normal to second order, and in usual settings these are
second-order equivalent to $r^*_\psi$. Moreover, it follows from
Severini \cite{seve2} that for likelihood-like objects that
are not true likelihoods, the first two Bartlett identities are enough
to validate the
$\mathrm{NP}$
adjustment referred to above.

The proposed methods apply directly to stratified Cox regression, by
adapting within each stratum to the observed censoring configuration.
An \texttt{R} package and \texttt{STATA} routine are available at
\url{http://www.science.oregonstate.edu/\textasciitilde piercedo/}.
For either package, the required user effort is essentially the same as
for ordinary Cox regression.

\section*{Acknowledgements}
The first author's work was partly supported by the Oregon Clinical and
Translational Research Institute at his university, funded by the
National Institutes of Health. Both authors were partly supported by
grants to the second author's university from the
\emph{Ministero dell'Universit\`{a} e della Ricerca}.
We are grateful to an Associate Editor for supporting this publication,
and to two anonymous referees for valuable suggestions. Dawn Peters was
very helpful in carefully reading the draft and offering comments that
improved the presentation.

%suskaldyti doi

% imsref loaded by jurgita.kaciuliene, 2014-02-27 10:57:19

\printhistory

\end{document}